\documentclass{amsproc}
\usepackage{epsfig}
\usepackage{graphicx}

\newtheorem{theorem}{Theorem}[section]

\theoremstyle{definition}

\newtheorem{example}[theorem]{Example}

\newtheorem{conjecture}[theorem]{Conjecture}

\theoremstyle{remark}

\newcommand{\figref}[1]{Fig.~\protect\ref{#1}}

\numberwithin{equation}{section}

\newcommand{\CL}{{\mathcal L}}

\newcommand{\CO}{{\mathcal O}}

\def\IZ{{\mathbb Z}}

\def\IC{{\mathbb C}}
\def\IP{{\mathbb P}}

\newcommand{\tr}{{\rm Tr}}
\newcommand{\re}{{\rm e}}

\newcommand{\rd}{{\rm d}}

\newcommand{\be}{\begin{equation}}
\newcommand{\ee}{\end{equation}}
\newcommand{\ba}{\begin{aligned}}
\newcommand{\ea}{\end{aligned}}
\newcommand{\ben}{\begin{eqnarray}\displaystyle}
\newcommand{\een}{\end{eqnarray}}

\newcommand{\sectiono}[1]{\section{#1}\setcounter{equation}{0}}

\newdimen\tableauside\tableauside=1.0ex
\newdimen\tableaurule\tableaurule=0.4pt
\newdimen\tableaustep
\def\phantomhrule#1{\hbox{\vbox to0pt{\hrule height\tableaurule width#1\vss}}}
\def\phantomvrule#1{\vbox{\hbox to0pt{\vrule width\tableaurule height#1\hss}}}
\def\sqr{\vbox{%
  \phantomhrule\tableaustep
  \hbox{\phantomvrule\tableaustep\kern\tableaustep\phantomvrule\tableaustep}%
  \hbox{\vbox{\phantomhrule\tableauside}\kern-\tableaurule}}}
\def\squares#1{\hbox{\count0=#1\noindent\loop\sqr
  \advance\count0 by-1 \ifnum\count0>0\repeat}}
\def\tableau#1{\vcenter{\offinterlineskip
  \tableaustep=\tableauside\advance\tableaustep by-\tableaurule
  \kern\normallineskip\hbox
    {\kern\normallineskip\vbox
      {\gettableau#1 0 }%
     \kern\normallineskip\kern\tableaurule}%
  \kern\normallineskip\kern\tableaurule}}
\def\gettableau#1{\ifnum#1=0\let\next=\null\else
\squares{#1}\let\next=\gettableau\fi\next}

\def\IE{\mathbb{E}}

\begin{document}

\title{Hurwitz numbers, matrix models and enumerative geometry}
\author{Vincent Bouchard}
\address{Jefferson Physical Laboratory, Harvard University,
17 Oxford St., Cambridge, MA 02138, USA}
\email{bouchard@physics.harvard.edu}
\author{Marcos Mari\~no}
\address{Department of Physics, Theory Division, CERN,
Geneva 23, CH-1211 Switzerland}
\email{marcos@mail.cern.ch}
\thanks{We would like to thank R. Pandharipande and R. Vakil for useful discussions, and A. Klemm and
S. Pasquetti for collaboration in a related work. Preprint number: CERN-PH-TH/2007-152.}

\subjclass{Primary 14N35; Secondary 14N10, 81T30}


\keywords{String theory, Hurwitz theory, matrix models}

\begin{abstract}
We propose a new, conjectural recursion solution for Hurwitz numbers at all genera. This conjecture is based
on recent progress in solving type B topological string theory on the mirrors of toric Calabi--Yau manifolds,
which we briefly review to provide some background for our conjecture. We show in particular
how this B-model solution, combined with mirror symmetry for the one-leg, framed topological vertex, leads to a
recursion relation for Hodge integrals with three Hodge class insertions. Our conjecture in Hurwitz theory follows from this recursion for the framed vertex in the limit of infinite framing.
\end{abstract}

\maketitle



\section{Introduction}

This article has two parts that can be read independently, although they are logically related to each other. The first part presents a new conjecture in Hurwitz theory. Hurwitz numbers can be regarded as the simplest enumerative quantities in algebraic
geometry, since they count covers of Riemann surfaces by Riemann surfaces. Recently, fascinating connections have emerged between
this classical problem and the most modern incarnations of enumerative geometry. Hurwitz numbers turn out to be closely related to
topological gravity in two dimensions \cite{gjv,opone}, to Hodge integrals \cite{elsv,fp,gvak}, and to the Gromov--Witten theory of $\IP^1$ \cite{optwo}.

An interesting problem in Hurwitz theory has been to find a complete set of recursion relations
solving the theory. Such relations are known to exist, and they can take many forms; see \cite{vakilrecursion,gj,gjv, gjva,fp,todapan} for significant examples. In this paper, we provide a new conjectural recursion solution of Hurwitz theory, different from the ones that have been presented in the literature, which
determines completely the simple Hurwitz numbers, counting covers of $\IP^1$ with arbitrary genus and with arbitrary
ramification at one point. Our result should be regarded as a conjecture, since we do not provide a proof.
We try however to state it in a precise way
and we give what we consider to be convincing evidence.

This conjecture in Hurwitz theory is in fact a spinoff of some recent progress
in topological string theory/Gromov--Witten theory on toric manifolds. In \cite{mm,bkmp}, a new formalism for the type B topological string
on mirrors of toric Calabi--Yau threefolds was proposed, inspired by the results of \cite{eo}. The recursion relations for
Hurwitz numbers that we obtain are in fact modeled on the
recursion procedure discovered in \cite{eo} in the context of matrix models. This procedure
has already found applications in two-dimensional topological gravity, and it has been
shown in \cite{eomir,emir} that Mirzakhani's recursion relations \cite{mir}
can be regarded as particular cases of the general formalism of \cite{eo}.

The second part
of the paper is devoted to a brief survey of this progress in matrix model theory and topological string theory, intended for a
mathematical audience. In particular, we review the proposal of \cite{mm,bkmp} that the one-leg, framed topological vertex obeys a similar set of
recursion relations by using the mirror geometry first found in \cite{akv}. This result, together with the
conjectural formula in \cite{mv} (later proved in \cite{llz,ophodge}) leads to a recursion relation
for Hodge integrals with three Hodge class insertions. It is known that Hurwitz numbers can be obtained as a limit of
the one-leg framed vertex in the case of infinite framing, and this finally leads to the
conjectural recursion
relation for Hurwitz numbers.

The organization of this paper is as follows. In section 2 we present our conjectural recursion relation in Hurwitz theory. In section 3 we
explain how this recursion arises in the context of mirror symmetry and the reformulation of the B-model proposed in \cite{mm,bkmp}.

\section{A new conjecture in Hurwitz theory}

\subsection{Hurwitz numbers}

In this paper we will be concerned with Hurwitz numbers counting the number of covers of $\IP^1$ by a Riemann surface of genus $g$ and
with arbitrary ramification at one point. The ramification type of the covering map at the special point is labeled by a partition 
$\mu$ of length $h:=\ell(\mu)$:
\be
\mu=(\mu_1 \geq \cdots \geq \mu_{\ell(\mu)} > 0).
\ee
We will use the standard notation
\be
|\mu|= \sum_{i=1}^{\ell(\mu)} \mu_i, \qquad z_\mu =  | {\rm Aut}(\mu) | \prod_{i=1}^{\ell(\mu)} \mu_i.
\ee
A partition $\mu$ can be regarded as a conjugacy class of the group of permutations of $|\mu|$ elements, $S_{|\mu|}$.
Given a representation $R$ of $S_{|\mu|}$, we define,
\be
\label{fr}
f_R(\mu)= {|\mu|!\over z_\mu} {\chi_R(\mu) \over{\rm dim}\, R},
\ee
where ${\rm dim}\, R$ is the dimension of $R$ and $\chi_R(\mu)$ is the character of the conjugacy class $\mu$ in the representation $R$. Representations $R$ of the symmetric group
will be identified by Young tableaux, labeled by the number of boxes $l_i$ in each row. The total number
of boxes of $R$ will be denoted by
\be
|R|=\sum_i l_i.
\ee
We will set $f_R(\mu)=0$ if $|R|\not=|\mu|$. We also define the following quantity associated to a tableau $R$,
\be
\label{kapar}
\kappa_R=\sum_i  l_i(l_i-2i +1).
\ee

A classical result expresses the Hurwitz numbers $H_{g,\mu}^{\bullet}$ counting disconnected covers
in terms of the representation theory of the
symmetric group $S_{|\mu|}$,
\be
H_{g,\mu}^{\bullet}=\sum_{R}\biggl( {{\rm dim}\, R \over |\mu|!}\biggr)^2 f_R(\mu) (\kappa_R/2)^{2g-2+h +|\mu|},
\ee
where only representations with $|R|=|\mu|$ contribute to the sum.
We also introduce a generating functional for Hurwitz numbers,
\be
\label{genhgf}
Z(v)=1+ \sum_{\mu} \sum_{g\ge 0} g_s^{2g-2+\ell(\mu)} { H_{g,\mu}^{\bullet} \over (2g-2+\ell(\mu) + |\mu|)!} p_{\mu}(v),
\ee
where $v=\{v_i\}_{i\ge 1}$ is a set of infinite formal variables and $p_{\mu}(v)$ are the symmetric polynomials
\be
p_{\mu}(v)=p_{\mu_1}(v) \cdots p_{\mu_h}(v), \qquad p_r=\sum_i v_i^r,
\ee
i.e. $p_r$ is the $r$th power sum (see \cite{macdonald}).
We will call $Z(v)$ the {\it partition function} of Hurwitz theory. It can be also written as
\be
Z(v)= \sum_R g_s^{-|R|}  \biggl( {{\rm dim} \, R \over |R|!}\biggr) \re^{-g_s \kappa_R/2} s_R (v),
\ee
where $s_R(v)$
is the Schur polynomial.
One can extract the formal logarithm of this, to obtain the {\it free energy} of Hurwitz theory,
\be
F(v)= \sum_{\mu} \sum_{g\ge 0} g_s^{2g-2+\ell(\mu)} { H_{g,\mu} \over (2g-2+\ell(\mu) + |\mu|)!} p_{\mu}(v),
\ee
which defines the Hurwitz numbers $H_{g,\mu}$ counting connected covers.

We will find it very convenient to collect the Hurwitz numbers at
fixed $g$ and $h=\ell(\mu)$ in the following functionals
\be \label{genh} H_g(x_1, \dots, x_h)=\sum_{\mu|\ell(\mu)=h} {z_\mu
\over (2g-2+  \ell(\mu)+|\mu|)!} H_{g,\mu} m_{\mu}(x), \ee
where $m_{\mu}(x)$ are the monomial symmetric functions in the $x_i$, $i=1, \cdots, h$, divided by an overall factor $x_1 \cdots x_h$,
\be
\label{msp}
m_{\mu}(x) = {1 \over |{\rm Aut} (\mu) |} \sum_{\sigma \in S_h} x_{\sigma(i)}^{\mu_i-1}.
\ee
We have that
\be
H_g(x_1, \cdots, x_h)={h!\over (2g-2+2h)!}H_{g,{\bf 1^h}}^{\IP^1} +\CO(x_1, \cdots, x_h).
\ee

We can use the ELSV formula to write the generating functionals (\ref{genh}) in a nice form (see \cite{gjv} for a similar rearrangement). Let $\overline{M}_{g,h}$ be the Deligne-Mumford moduli space of genus $g$ curves with $h$ marked points, $\IE$ the Hodge bundle, and $\lambda_i = c_i(\IE)$ its Chern classes. Define
\be
\Lambda^\vee_g (t) = \sum_{i=0}^g (-1)^i \lambda_i t^{g-i}.
\ee
The ELSV formula states that
\be
\label{elsvfor}
H_{g,\mu} =
{ (2g-2+ \ell(\mu)+|\mu|)! \over z_\mu} \prod_{i=1}^{\ell(\mu)} {\mu_i^{\mu_i+1}\over \mu_i!} \int_{\overline{M}_{g,\ell(\mu)}}
{\Lambda^\vee_g(1)
\over \prod_{i=1}^{\ell(\mu)} (1-\mu_i \psi_i)}.
\ee
The $\psi_i$ classes of two-dimensional topological gravity are
defined as follows. Let $\CL_i$ be the line bundle over
$\overline{M}_{g,h}$ whose fiber at a point $(C; x_1, \dots, x_h)
\in \overline{M}_{g,h}$ is the cotangent space to $C$ at $x_i$.
Then $\psi_i := c_1(\CL_i)$.

We now introduce the variable $y$ through
\be
\label{xyrel}
x=y\re^{-y},
\ee
which defines the so-called tree function $T(x)$
\be
y(x)=T(x) =-W(-x)
\ee
where $W(x)$ is the Lambert $W$ function (see \cite{lambert} for a very useful summary of its properties). Near $(x,y)=(0,0)$
one has the series expansion
\be
\label{trees}
y(x)=\sum_{\mu=1}^{\infty} {\mu^{\mu-1} \over \mu!} x^{\mu}.
\ee
We now define the one-forms
\be
\zeta_n(y)=\rd y {1-y \over y} \Bigl( {y \over 1-y} {\rd \over \rd y}\Bigr)^{n+2} y, \qquad n\ge 0.
\ee
For example,
\be
\zeta_0(y) ={\rd y \over (1-y)^2}, \quad \zeta_1(y) = {1+ 2y \over (1-y)^4} \rd y.
\ee
Since
\be
{y\over1- y}  {\rd \over \rd y}=x {\rd \over \rd x}, \quad {\rd x\over x}= {1-y \over y} \rd y,
\ee
we have the following expansion in powers of $x$
\be
\label{zeta1}
\zeta_n(x)  ={\rd x  \over x} \Bigl( x{\rd \over \rd x}\Bigr )^{n+2} y(x)=\sum_{\mu=1}^{\infty} {\mu^{\mu+1+n} \over \mu!} x^{\mu-1} \rd x.
\ee
It then follows that
\be
\ba
\label{elsvgen}
& H_{g}(x_1, \dots, x_h) \rd x_1 \cdots \rd x_h  \\
&= \sum_{n_i\ge0} \sum_{\mu_i \geq 1} \left( \prod_{i=1}^h {\mu_i^{\mu_i+1+n_i }\over \mu_i!}x_i^{\mu_i-1} \rd x_i\right)  \int_{\overline{M}_{g,h}}
\Lambda_g^\vee(1)  \prod_{i=1}^h \psi^{n_i}_i\\ &=
\sum_{n_i=0}^{3g-3+h} \langle \tau_{n_1} \cdots \tau_{n_h} \Lambda_g^\vee (1) \rangle
\prod_{i=1}^h \zeta_{n_i}(y_i),
\ea
\ee
where we introduced the Hodge integrals
\be
\langle \tau_{n_1} \cdots \tau_{n_h} \Lambda_g^\vee (1)\rangle = \int_{\overline{M}_{g,h}}
\Lambda_g^\vee (1) \prod_{i=1}^h \psi^{n_i}_i.
\ee

\subsection{A recursion relation for Hurwitz theory}

Our conjecture will describe a recursion solution for the generating functionals of Hurwitz numbers \eqref{genh}. This recursion relation was first proposed by Eynard and Orantin in a
different context, namely for correlation functions of matrix models \cite{eo}. Then, in \cite{mm,bkmp} it was conjectured that a suitable generalization of this relation also describes B-model topological string theory on the mirrors of toric Calabi-Yau threefolds. In the present paper we conjecture that this relation also generates recursively the generating functionals of Hurwitz numbers.

We will review these developments in section 3, and motivate our
conjecture for Hurwitz numbers in section 3.4 by using the
relation between Gromov-Witten theory and Hurwitz theory. We will
see that the recursion relation in Hurwitz theory is in fact a
particular case of a similar recursion relation describing
Gromov-Witten potentials (or topological string amplitudes) on
toric Calabi--Yau threefolds. But, for the moment, let us explain
the recursion in the context of Hurwitz theory and state our
conjecture.

Generically, the recursion relations of Eynard and Orantin provide a formalism to generate an infinite sequence of meromorphic differentials and symplectic invariants associated to a curve. Our application to Hurwitz theory consists in a particular case of this construction, where we consider the curve introduced earlier,
\be
C: \{ x = y \re^{-y} \},
\ee
which defines the tree function $T(x)$
\be
y(x) = T(x).
\ee

To formulate the recursion relations we need a few geometric ingredients associated to the curve $C$.
First, the $x$-projection of $C$ has one ramification point, which we will denote by $\nu$, at
\be
(x(\nu), y(\nu) ) = (\re^{-1}, 1).
\ee
Near the ramification point $\nu$ there are two points $q, \bar q
\in C$ with the same $x$-projection, that is, $x(q) = x(\bar q)$.
Write
\be
y(q) = 1 + z, \qquad y(\bar q) = 1 + S(z),
\ee
with
\be
S(z) = -z + \CO(z^2).
\ee
Note that in these variables, the ramification point is at $z=0$. By definition of $q$ and $\bar q$, we have
\be
(1+z) \re^{-z} = (1+S(z) ) \re^{-S(z)},
\ee
which we can solve to get $S(z)$ as a power series in $z$:
\be
S(z) = -z + \frac{2\,z^2}{3} - \frac{4\,z^3}{9} + \frac{44\,z^4}{135} - \frac{104\,z^5}{405} + \frac{40\,z^6}{189} - \frac{7648\,z^7}{42525}+\CO(z^8).
\ee
This series arises in the study of random graphs; see for example
\cite{knuth}.

The first ingredient is the Bergman kernel of $C$, which is the unique
meromorphic differential $B(x_1,x_2)$ on $C$ with a double pole at $x_1=x_2$. Since $C$ has genus $0$, the Bergman kernel is given, in terms of the local coordinate $y$, by
\be
\label{Bergman}
B(y_1,y_2)= {\rd y_1 \rd y_2 \over (y_1-y_2)^2 },
\ee
where $y_i$ is implicitly related to $x_i$ by the tree function,
that is, $y_i = y(x_i) = T(x_i)$. The next ingredient is built up
from this, and it is the one-form
\be
\rd E_q(x)= {1 \over 2} \int_{q}^{\bar q} B(\xi,x),
\ee
which is well-defined locally near the ramification point $\nu$. For $C$, we have, in the coordinates $z$ and $y$,
\be
\label{deq}
\rd E_z(y)= {1 \over 2} \rd y \Bigl[{1\over y -1 - z}-{1\over y -1 - S(z)}\Bigr].
\ee
Finally, we need the one-form
\be
\omega(q)=(\log y(q)- \log y(\bar q)) {\rd x (q) \over x(q) },
\ee
which becomes, in terms of $z$,
\be
\ba
\omega(z) &= \Big(\log (1+S(z)) - \log (1+z) \Big) {z \rd z \over 1+z } \\
&= (S(z) - z ) {z \rd z \over 1+z }.
\ea
\ee

\begin{figure}
\leavevmode
 \begin{center}
\epsfxsize=5in
\epsfysize=1.3in
\epsfbox{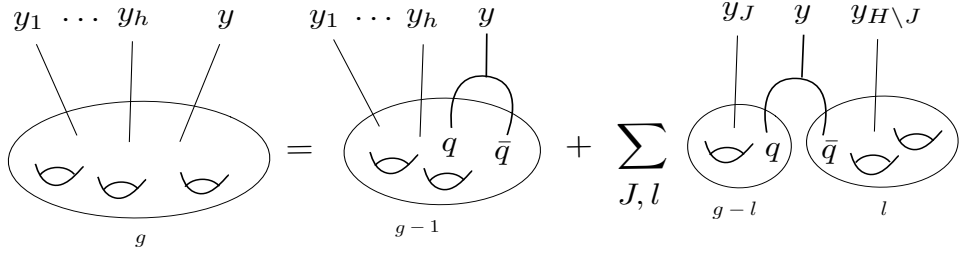}
\end{center}
\caption{A graphic representation of the recursion relation (\ref{recursion}).}
\label{recfig}
\end{figure}

With these ingredients, we define recursively a set of meromorphic
differentials $W_g(y_1, \dots, y_h)$, with $g\ge 0$, $h\ge 1$,
using residue calculus at the ramification point $z=0$ of the
$x$-projection of the curve $C$. We first define the starting
point of the recursion as
\be
W_0(y) = 0, \qquad W_0(y_1,y_2) = B(y_1,y_2).
\ee
The meromorphic differentials are then defined by
\be
\ba
W_g(y, y_1 \dots, y_h) &= {\rm Res}_{z=0}{\rd E_{z}(y) \over \omega(z)} \Bigl[ W_{g-1} (1+z, 1+S(z), y_1, \dots, y_h)\\
&\quad +\sum_{l=0}^g \sum_{J \subset H} W_{g-l}(1+z, y_J) W_l(1+S(z), y_{H\backslash J})\Bigr],
\label{recursion}
\ea
\ee
where $H={1, \dots, h}$, and given any subset $J=\{i_1, \dots,
i_j\}\subset H$ we set
\be y_J=\{y_{i_1}, \dots, y_{i_j}\}. \ee
The recursion relation \eqref{recursion} can be represented graphically as in \figref{recfig}.

We are now ready to state our main conjecture:

\begin{conjecture}

The generating functionals of Hurwitz theory in the differential
form \eqref{elsvgen} are explicitly given as follows (we use the
notation above). For $g=0, h=1$ and $g=0, h=2$ one has
\be
H_0(x)\rd x = y(x){\rd x \over x} , \qquad H_0(x_1, x_2) \rd x_1 \rd x_2 =B(y_1, y_2) -{\rd x_1 \rd x_2 \over (x_1-x_2)^2},
\ee
where $y(x) = T(x)$ and the Bergman kernel $B(y_1,y_2)$ was defined in \eqref{Bergman}.
For the remaining cases, one has that
\be H_g(x_1, \dots, x_h)\rd x_1 \cdots \rd x_h =W_g(y_1, \dots,
y_h) \ee
where $W_g(y_1, \dots, y_h)$ are defined by the recursion relation
(\ref{recursion}). \label{conj1}
\end{conjecture}

%
%

This conjecture will be further motivated in section 3.4, using the relations between Hurwitz theory, Gromov--Witten theory and B-model topological string theory. Let us now present some examples of the generating functionals.

We find, for example, that
\be
\ba
W_1(y)&=-{1\over 2} {\rm Res}_{z=0} \Biggl[ {1+z \over z} \Bigl({1\over y-1-z} -{1\over y-1-S(z)}\Bigr) {S'(z) \rd z
\over (z-S(z))^3}\Biggr] \rd y \\
&={4y -y^2\over 24 (1-y)^2} \rd y.
\ea
\ee
When expanded in terms of $x$, we recover known results on Hurwitz numbers.
Notice that it is more convenient to work with the $y$ variable, rather than the $x$ variable, and this is in
turn very useful in order to extract the Hodge integrals in \eqref{elsvfor} from the results, since one only has to write them in terms of the
one-forms $\zeta_n(y)$ and compare them to (\ref{elsvgen}). For instance,
\be
W_1(y)={1 \over 24}\bigl( -\zeta_0(y) +\zeta_1 (y)\bigr) .
\ee
Similarly, one can find by using the residue calculus,
\be
\ba
W_0(y_1, y_2, y_3)&=\prod_{i=1}^3   \zeta_0(y_i) , \\
W_0(y_1, \dots, y_4)&= \sum_{i=1}^4 \zeta_1(y_i) \prod_{j\not= i}
\zeta_0(y_j), \ea \ee
for $g=0$. For $g=1, h=2$ the result of the recursion gives
\be
\ba
W_1(y_1, y_2)&={1 \over 24}\Bigl( -\zeta_0(y_1) \zeta_1(y_2)+ \zeta_0(y_1) \zeta_2(y_2) +(y_1 \leftrightarrow y_2)
\\ & \qquad \qquad
+\zeta_1(y_1) \zeta_1(y_2) \Bigr) .
\ea
\ee
For $g=2$ and $h=1,2$, one obtains,
\be
\ba
W_2(y) &= {1 \over 5760} \Big( {7} \zeta_2(y)- 12 \zeta_3(y) +{5} \zeta_4(y) \Big),\\
W_2(y_1, y_2)&= {1 \over 5760} \Big({7} \zeta_0(y_1) \zeta_3(y_2) -{12}  \zeta_0(y_1) \zeta_4(y_2) +{5}  \zeta_0(y_1) \zeta_5(y_2) \\
&\qquad\qquad +{21}  \zeta_1(y_1) \zeta_2(y_2)  -{36} \zeta_1(y_1) \zeta_3(y_2) +{15} \zeta_1(y_1) \zeta_4(y_2) \\
&\qquad\qquad +{29} \zeta_2(y_1) \zeta_3(y_2) +(y_1 \leftrightarrow y_2)
 -{50} \zeta_2(y_1) \zeta_2(y_2) \Big) ,
 \ea
 \ee
 and finally for $g=3$, $h=1$ one has
 \be
W_3(y)= {1 \over 2903040} \Big( -{93} \zeta_4(y) + {205}\zeta_5(y) -{147} \zeta_6(y)
+ {35} \zeta_7(y) \Big).
\ee
By comparing these results with the expected form \eqref{elsvgen}, it is easy to check that the recursion
relation reproduces known results for the Hodge integrals involved in \eqref{elsvfor}. This implies in
particular that the above conjecture provides a recursion relation for the Hodge integrals appearing
in the ELSV formula (which include all the correlation functions of 2d gravity).

Some remarks concerning these results are in order:

\begin{enumerate}

\item In view of the results of \cite{eomir,emir}, our conjectural recursion formula is probably
the analogue in Hurwitz theory of Mirzakhani's recursion for the Weil--Petersson volume of moduli space \cite{mir}.
It is known that Mirzakhani's recursion is essentially equivalent to the KdV/Virasoro structure of 2d topological gravity \cite{mulase,mirtwo,liuzu}.
In the same way, we expect our recursion relation to be equivalent to the Toda structure governing Hurwitz theory \cite{todapan,okounkov}.

\item It may be worth mentioning that the conjectural recursion that we propose is different from the cut-and-join equation for simple Hurwitz numbers (see for example
Lemma 3.1 in \cite{gjva}). The latter is more easily formulated in the ``representation basis" (consisting of the symmetric polynomials $p_\mu(v)$), while our recursion is naturally formulated in the ``winding number basis" (consisting of the monomial symmetric functions $m_\mu(x)$); moreover, the initial conditions of the recursions are different. In fact, the cut-and-join equation should be regarded as a differential equation governing the generating functional of Hurwitz numbers, while our recursive conjecture should be regarded as an explicit solution of this equation. Nevertheless, it should be possible to derive one from the other. In particular, it would be very interesting to prove our recursion using the cut-and-join equation as a starting point. However, the relation between the two equations seems to involve highly non-trivial combinatorics, which appear to be nicely encoded in the function $S(z)$ and the residue calculus.

\end{enumerate}

\sectiono{Toric geometry, matrix models, and mirror symmetry}

To provide further motivation for the recursion conjecture \ref{conj1} in Hurwitz theory, we need to introduce a new framework for B-model topological string theory on mirrors of toric Calabi--Yau threefolds, inspired by matrix models. Let us start by reviewing known results in the matrix model realm.

\subsection{Matrix models}

The recursion relation \eqref{recursion} which we use with small modifications in Hurwitz theory was first
found in the context of matrix models. A matrix model for a Hermitian $N \times N$ matrix $M$ is defined by a partition function
\be \label{zmm} Z=\int \rd M \, \re^{-N \tr \, V(M)}, \ee
where $V(x)$ is a polynomial, or more generally a power series in $x$. The main goal of matrix model
theory is to compute \eqref{zmm} as well as the connected correlation functions
\be \label{wmm} W(z_1, \dots, z_h)=\Big\langle  {\rm Tr}\, {1\over
z_1 -M} \cdots  {\rm Tr}\, {1\over z_h -M}\Big\rangle^{(c)}. \ee
It turns out that these quantities have an asympotic expansion in inverse powers of $N$ (the rank of the matrix), with the structure
\be
\ba
\log Z&=\sum_{g=0}^{\infty} F_g N^{2-2g},\\
W(z_1, \dots, z_h)&=\sum_{g=0}^{\infty} N^{2-2g+h} W_g(z_1,
\dots, z_h). \ea \ee
We will call $F_g$ the genus $g$ closed amplitude, and $W_g(z_1,
\dots, z_h)$ the genus $g$ open amplitude with $h$ holes. The
reason for these names is their diagrammatic interpretation in
terms of double-line diagrams or fatgraphs (see for example
\cite{biz} for an exposition of the diagrammatics of matrix
models).

The $1/N$ expansion of matrix models has been studied for thirty
years now, starting with the seminal paper \cite{bipz}. The main
result that emerges is that the genus $g$ amplitudes can be
obtained from a single object, an algebraic curve \be \{ P(x,y)=0
\} \subset \IC^2 \label{spectral} \ee called the {\it classical
spectral curve} of the matrix model. This curve is closely related
to the generating function of disk amplitudes (also called the
resolvent):
\be
W_0(x)={1\over 2}(V'(x) + y(x)).
\ee
It turns out that the genus $g$ open amplitudes of the matrix model can be constructed recursively as follows (see \cite{eo} for a detailed exposition and a list of references).\footnote{In fact, the genus $g$ amplitudes constructed in this manner on any algebraic curve in $\IC^2$ provide symplectic invariants of the curve, whether the curve is the spectral curve of a matrix model or not. This was the main insight of \cite{eo}.}

We first need some necessary ingredients, which were introduced in section 2.2 for the curve relevant for Hurwitz theory. To start with, we notice that, generically, the $x$-projection of the
curve $C$ will have a set of {\it ramification points} $q_i \in C$ defined by the condition
\be
\label{rama}
\rd x (q_i)=0.
\ee
Near a ramification point\footnote{We consider only simple ramification points.} there are two points, $q$ and $\bar q$, with the same $x$-coordinate $x(q)=x(\bar q)$.
The first ingredient is the {\it Bergman kernel} of $C$, which is the unique
meromorphic differential on $C$ with a double pole at $p=q$, and normalized such that
\be
\oint_{A^I} B(p,q) =0,
\ee
where $(A^I, B_I)$ is a canonical basis of cycles on the Riemann surface. If $C$ has genus $0$ (which is the only case we will
consider in some detail in this paper), it is given explicitly by, in terms of the local coordinate $y$,
\be
B(p,q) = { \rd y(p) \rd y(q) \over (y(p)-y(q))^2 }.
\ee
The next ingredient is the one-form
\be
\rd E_q(p)= {1 \over 2} \int_{q}^{\bar q} B(\xi,p),
\ee
which is well-defined locally near a ramification point $q_i$. In genus $0$, we have
\be
\label{deqgen}
\rd E_q(p)= {1 \over 2} \rd y(p) \Bigl[{1\over y(p) -y(q)}-{1\over y(p) -y(\bar q)}\Bigr].
\ee
Finally, we need the one-form
\be
\omega(q)=(y(q)-y(\bar q)) \rd x (q).
\ee
With these ingredients, \cite{eo} define recursively a set of
meromorphic differentials $W_g(p_1, \dots, p_h)$, with $g\ge 0$,
$h\ge 1$, as follows. We first define the starting point of the
recursion as
\be
W_0(p) = 0, \qquad W_0(p_1,p_2) = B(p_1,p_2).
\ee
The meromorphic differentials are then defined by
\be
\ba
W_g(p, p_1 \dots, p_h) &= \sum_{q_i} {\rm Res}_{q=q_i}{\rd E_{q}(p) \over \omega(q)} \Bigl[ W_{g-1} (q, \bar q, p_1, \dots, p_h)\\
&\quad +\sum_{l=0}^g \sum_{J \subset H} W_{g-l}(q, p_J) W_l(\bar q, p_{H\backslash J})\Bigr],
\label{recursion2}
\ea
\ee
where the sum is over all ramification points defined by (\ref{rama}), and the notation is as in \eqref{recursion}.

Note that it is also possible to calculate the closed amplitudes $F_g$ by the following formula,
\be
F_g=\sum_i {\rm Res}_{p=q_i} \Phi(q) W_g(q),
\ee
where $\Phi(q)$ is any primitive of $y(q)$.

\subsection{Topological strings on toric Calabi--Yau threefolds}

Topological string theory on a Calabi--Yau manifold $X$ comes in two varieties, called the A- and the B-models. In the A-model,
the genus $g$ closed amplitudes $F_g$ are generating functionals for
Gromov--Witten invariants $N_{g,Q}$ ``counting" holomorphic maps from genus $g$ Riemann surfaces $\Sigma_g$ to $X$
\be
F_g(t) =\sum_{Q\in H_2(X)} N_{g,Q} \re^{-Q\cdot t}
\ee
where $t$ are complexified K\"ahler parameters.
In the B-model we study variation of complex structures on a mirror Calabi--Yau manifold $\widetilde X$.
The basic buiding blocks are the period integrals of the
holomorphic 3-form $\Omega$ on $\widetilde X$
\be
\label{compactperiods}
t^I =\oint_{A^I} \Omega, \qquad {\partial F_0 \over \partial t^I} = \oint_{B_I} \Omega
\ee
where $A^I, B_I$, $I=1, \dots, b_3$ is a symplectic basis of
$H_3(\widetilde X)$. This defines the genus zero closed amplitude
$F_0(t)$, where $t$ now parameterizes complex structures. The
genus one amplitude $F_1(t)$ can be defined in terms of the
Ray--Singer torsion of $\widetilde{X}$. The higher genus $F_g(t)$
are not well-defined mathematically, although they can be
partially constructed from the holomorphic anomaly equations of
\cite{bcov}.

\begin{figure}
\leavevmode
 \begin{center}
\epsfxsize=3.75in
\epsfysize=1.5in
\epsfbox{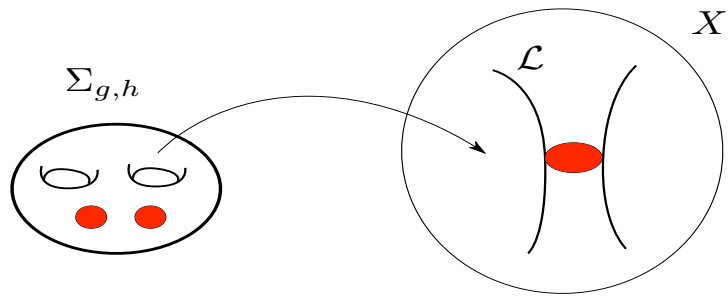}
\end{center}
\caption{In the presence of topological D-branes, we can consider maps from Riemann surfaces with boundaries.}
\label{dbranes}
\end{figure}

We can introduce an open sector in topological string theory by considering topological D-branes. Again, the description is easier to make in the
A-model, where the topological D-branes must wrap a Lagrangian submanifold $\CL \subset X$. In this way we can also consider
maps from Riemann surfaces
with genus $g$ and $h$ holes $\Sigma_{g,h}$ to $X$ with boundary conditions set by $\CL$, as shown in \figref{dbranes}.

Let us now assume that $b_1(\CL)=1$. The counting of maps is
encoded in a topological open string amplitude $F_{g,h}(x_i)$,
\be \label{fghgen} F_{g,h}(x_1, \dots, x_h)=\sum_{w_i} F_{g,w} (t)
x_1^{w_1} \cdots x_h^{w_h}, \ee
where the $w_i$ are the winding numbers of the holes of $\Sigma_{g,h}$ around the nontrivial cycle in $\CL$, and
\be
F_{g,w}(t) =\sum_{Q\in H_2(X)} N_{g,w,Q} \re^{-Q\cdot t}
\ee
counts holomorphic maps with boundary conditions set by $\CL$ in the topological sector specified by
$Q,w,g$. The quantities $N_{g,w,Q}$ are called open Gromov--Witten invariants; in general, they are not
well--defined mathematically. In the toric case, however, they can be fully formalized in terms of relative Gromov--Witten theory \cite{lllz}.
The variables $x_i$ in \eqref{fghgen} are called {\it open moduli}.

By mirror symmetry, there should be similar amplitudes in the B-model on $\widetilde{X}$, for appropriate objects which are mirror to
Lagrangian submanifolds (these objects are, in general, coherent sheaves). In general situations it is difficult to define these amplitudes --- although
in the case of rigid Lagrangian submanifolds in compact Calabi--Yau manifolds one can write down generalizations of the BCOV equations \cite{walcher}.
As we will see in a moment, in the toric case these open amplitudes can be defined by a recursion similar to (\ref{recursion2}).

\begin{figure}
\leavevmode
 \begin{center}
\epsfxsize=4.5in
\epsfysize=2in
\epsfbox{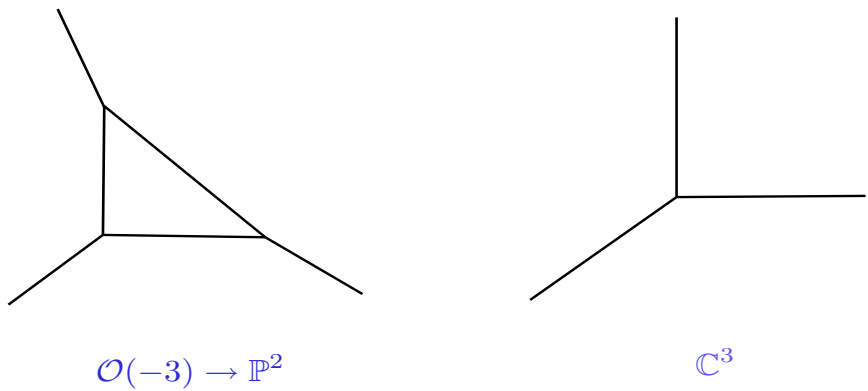}
\end{center}
\caption{Two simple examples of toric manifolds.}
\label{toricman}
\end{figure}

Let us then focus on the case in which $X$ is toric. A toric Calabi--Yau is necessarily noncompact, and can
be represented by a toric diagram which encodes its fibration structure (see for example \cite{mmrev} for details). In \figref{toricman} we
show some simple examples of toric diagrams: the total space of the bundle $\CO(-3) \rightarrow \IP^2$, also known as local $\IP^2$, and the even simpler case of
$\IC^3$.

In toric geometries there is an important class of D-branes wrapping Lagrangian submanifolds
with the topology $\IC\times S^1$, which we will call for short {\it toric branes}.
These Lagrangian submanifolds are generalizations to the toric realm of the
Lagrangian submanifolds of $\IC^3$ considered by Harvey and Lawson in \cite{hl}.
In the physics literature they were first studied in \cite{av} and further explored in \cite{akv}.
In the toric diagram these Lagrangian submanifolds
simply project to points in the edges, and one can study the open string amplitudes with the Lagrangian
boundary conditions that they define.

The A-model open and closed topological string amplitudes on toric Calabi--Yau threefolds can be computed with the so-called
topological vertex \cite{akmv}. In this formalism, the total closed partition function
\be
Z(t,g_s)=\exp \Bigl[ \sum_{g=0}^{\infty} F_g(t) g_s^{2g-2}\Bigr]
\ee
is obtained as a sum over partitions or Young tableaux,
\be
\label{ztv}
Z(t,g_s)=\sum_{R_i} C(R_i, g_s) \re^{-\sum_i \ell(R_i) t}.
\ee
In this expression, $i$ runs over the edges of the toric diagram, $\ell(R)$ is the length of the partition $R$, and $C(R_i,g_s)$ is a quantity that
depends on the representations, the variable $g_s$, and the shape of the toric diagram, and which can be computed in a very precise way by using
the theory of the topological vertex. The expression for the partition function in \eqref{ztv} is exact in $g_s$ (therefore it contains the all-genus information of the topological string), but it is perturbative in $\re^{-t}$. The expansion is around the large radius limit $t\rightarrow \infty$ of the K\"ahler
moduli space. For the open amplitudes, the result is perturbative in both $\re^{-t}$ and the open moduli $z_i$.

Therefore, the vertex formalism is only appropriate for some particular point in the open/closed moduli space (i.e. large radius for $t$, and the point $z_i=0$).
It is not well-suited for studying other points in the moduli space like the conifold and orbifold points. This is where the B-model becomes useful, since
it is nonperturbative in the K\"ahler parameter. Let us now review what is known about the B-model on the mirrors of toric manifolds.

The Calabi--Yau mirrors of toric geometries have been known for some time. Generically, they are given by conic fibrations over $\IC^* \times \IC^*$, where the conic fiber degenerates to two lines over an algebraic curve (see for example \cite{kkv,ckyz}, and \cite{hv} for a physics derivation). The geometry is essentially captured by the algebraic curve in $\IC^* \times \IC^*$, which is usually called the {\it mirror curve}, and is given by a polynomial
\be
\Sigma: \{H(x,y)=0\} \subset \IC^* \times \IC^*.
\ee
It is crucial to note that, in contrast with the spectral curve \eqref{spectral} of matrix models, the mirror curve is a Riemann surface embedded in $\IC^* \times \IC^*$ rather than $\IC^2$.

\begin{figure}
\leavevmode
 \begin{center}
\epsfxsize=4.5in
\epsfysize=2in
\epsfbox{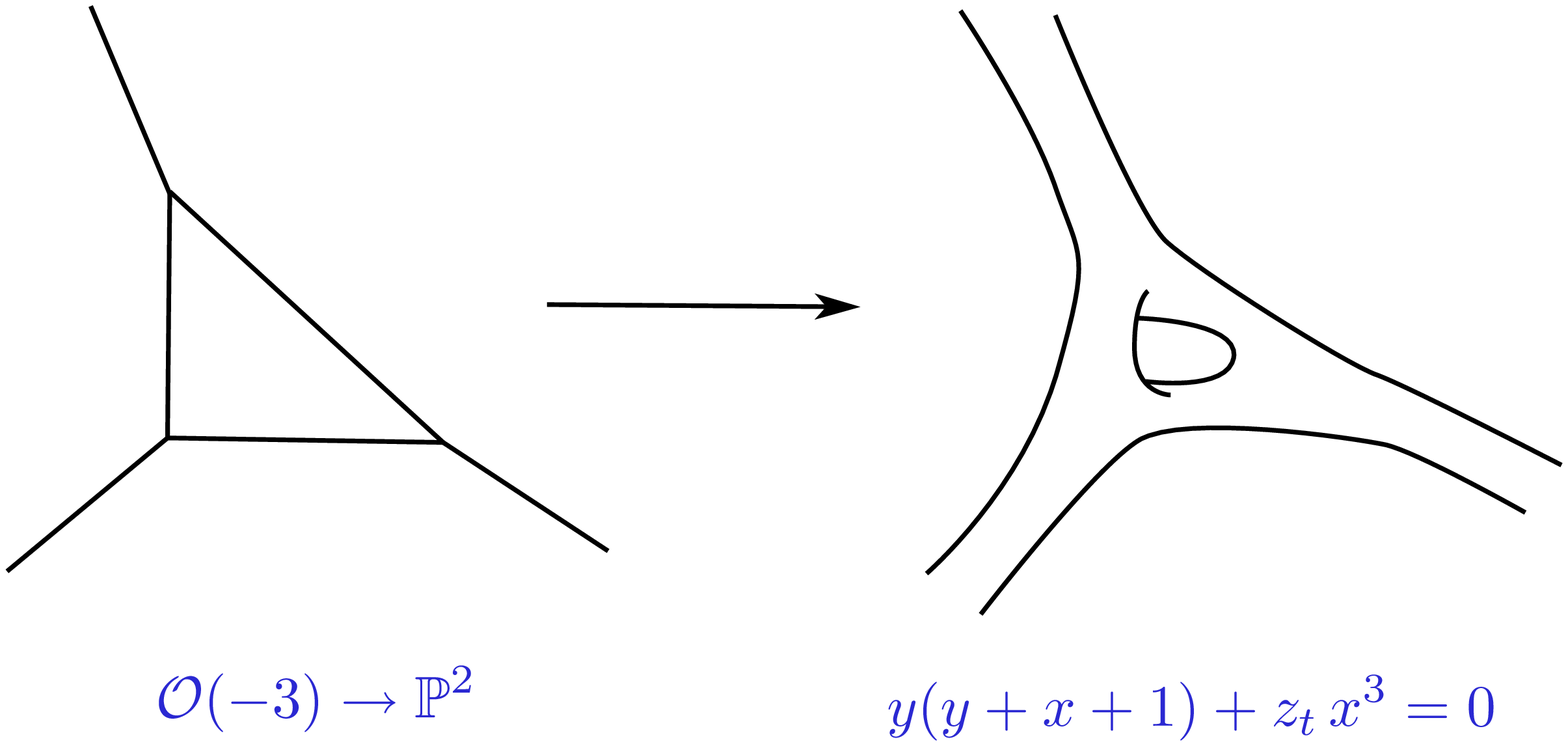}
\end{center}
\caption{Local $\IP^2$ and its mirror.}
\label{ptwomirror}
\end{figure}

The mirror curve of a toric Calabi--Yau threefold $X$ can be visualized as a thickening of the
toric diagram of $X$. The example of local $\IP^2$ and its mirror is shown in figure \figref{ptwomirror}; in this case, the mirror curve is a genus $1$ Riemann surface with three punctures.

The closed sector of the B-model on these geometries can be
analyzed by using techniques similar to those developed for the
compact case. In particular, the prepotential is determined by the
equations
\be
\label{logcompactperiods}
t^I =\oint_{A^I} \log y {\rd x \over x} , \qquad {\partial F_0 \over \partial t^I} = \oint_{B_I} \log y {\rd x \over x}
\ee
where $A^I, B_I$ is a symplectic basis of the mirror curve. Higher
genus amplitudes $F_g$ can be studied with a toric version of the
holomorphic anomaly equations; see for example \cite{kz}.

What about open string amplitudes? It turns out that the toric branes become, under mirror symmetry, just points on the mirror curve \cite{av}, and the
open string moduli appearing in (\ref{fghgen}) are local coordinates for the curve. The choice of local coordinate --- which is equivalent to a choice of projection $\Sigma \to \IC^*$ --- is related to the edge on which the brane is
located in the toric diagram of the A-model.

More precisely, as a curve in $\IC^* \times \IC^*$, the mirror curve $\Sigma$ has reparameterization group $SL(2,\IZ)$ acting as
\be
(x,y) \rightarrow (\tilde x, \tilde y)=(x^a y^b,x^c y^d),\qquad \begin{pmatrix}a & b\\
c&d \end{pmatrix} \in SL(2,\IZ).
\ee
Changing the parameterization of the curve corresponds to choosing different local coordinates (or open string moduli), hence moving the position of the toric brane on the toric diagram on the mirror side \cite{av,akv}. In fact, not all reparameterizations lead to different geometric branes; fixing the position of the brane only fixes the reparameterization of the curve up to a one-parameter subgroup of $SL(2,\IZ)$. This extra integer corresponds to the {\it framing} of the brane, which is an ambiguity in the computation of open string amplitudes. The framing transformations correspond to reparameterizations of the form
\be
(x,y)\mapsto (\tilde x, \tilde y) =  (xy^f,y), \qquad f \in \IZ.
\label{framing}
\ee

In \cite{av} it was further shown that the disk amplitude of a toric brane can be computed in terms of the defining equation of the mirror curve
\be
\label{disk}
F_{0}(x)\rd x=\log \, y {\rd x \over x},
\ee
which is nothing but the one-form appearing in (\ref{compactperiods}).

\begin{figure}
\leavevmode
 \begin{center}
\epsfxsize=4.5in
\epsfysize=2in
\epsfbox{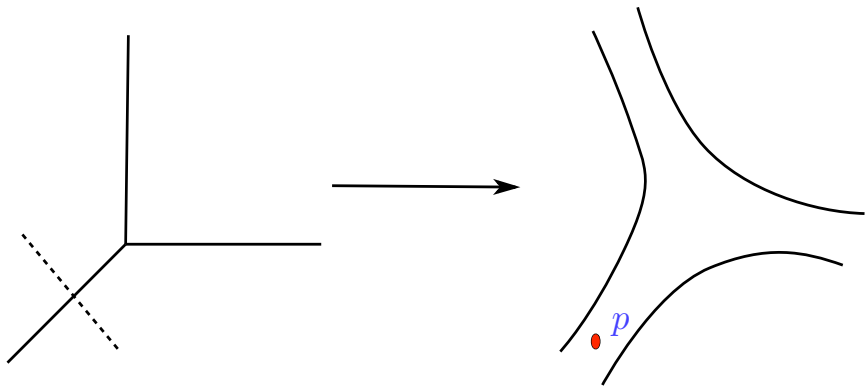}
\end{center}
\caption{$\IC^3$ with a brane in a leg, and its mirror.}
\label{cbrane}
\end{figure}

\begin{example}
\label{exc3}
A relevant example in the following will be $\IC^3$ with a brane in one leg. The mirror curve to $\IC^3$ reads
\be
1-y-x=0,
\label{mirrorc3}
\ee
which is $\IP^1$ with three punctures. The open modulus can be taken to be simply $x$, the local coordinate.
The D-brane configuration and its mirror are shown in \figref{cbrane}.

The framed brane in one leg of $\IC^3$, which we will call the {\it framed vertex}, can be obtained by applying the framing transformation \eqref{framing} to the mirror curve \eqref{mirrorc3}.
We get
\be
-\tilde y^{f+1} + \tilde y^f -\tilde x=0.
\label{framedc3}
\ee
The open modulus is now $\tilde x$, and the expression \eqref{disk} is still valid in the new parameterization of the curve given by $(\tilde x, \tilde y)$.
\end{example}

A natural question is how to compute $F_{g,h}(x_i)$ in the B-model, for arbitrary $g,h$. The main claim of \cite{mm} (backed by various
nontrivial examples), which was extensively developed and clarified in \cite{bkmp}, is that these amplitudes can be computed
by using a slight modification of the recursion relation \eqref{recursion2} first found in the context of matrix models. The classical spectral curve is now taken to be
the mirror curve. The only difference in the recursion relation comes from the fact that the mirror curve is embedded in $\IC^* \times \IC^*$ rather than $\IC^2$; consequently, the one-form $\omega(q)$ becomes
\be
\omega(q) = \big(\log y(q) - \log y(\bar q) \big) {\rd x(q) \over x(q)},
\label{oneform}
\ee
rather than $\omega(q) = (y(q) - y (\bar q) ) \rd x(q)$, which reflects the fact that the symplectic form on $\IC^* \times \IC^*$ is
\be
{\rd x \over x} \wedge {\rd y \over y}.
\ee

The physical reason for the claim in \cite{mm} is that the theory
of matrix models is just an example of a more general
construction, namely a chiral boson on a ``quantum" Riemann
surface, as was suggested in \cite{adkmv}. In that paper it was
also argued that the B-model on the mirror of a toric Calabi--Yau
threefold is another example of this theory, where the chiral
boson now lives  on the mirror curve. The fact that the invariants
defined in \cite{eo} depended only on the Riemann surface
indicated that they must be the amplitudes for the chiral boson
theory, which led to the claim of \cite{mm}. Note that after the
appearance of \cite{mm,bkmp}, this line of thought was further
refined in \cite{dvnew}, where a physical reinterpretation of the
recursion relation directly in terms of the chiral boson theory
was proposed.

From a more mathematical point of view, one can take the modified
recursion relation \eqref{recursion2} with the one-form
\eqref{oneform} to be the definition of the B-model on the mirrors
of toric geometries, and the claim of \cite{mm,bkmp} is then a
conjecture that these quantities correctly reproduce the open and
closed Gromov--Witten invariants of toric Calabi--Yau manifolds.
From this point of view, the recursive formalism can be understood
as a gluing procedure for open and closed topological string
amplitudes, where the basic building blocks are the disk and the
annulus amplitudes. A very particular case of this conjecture,
involving the disk amplitude of a certain class of branes in local
$\IP^2$, was already proved in \cite{gz}.

We will now consider in some more detail the simplest possible case of this formalism, namely the framed vertex introduced in example \ref{exc3}. As a spinoff we will derive our conjecture for the recursion relation of Hurwitz numbers.

\subsection{Framed vertex}

The A-model topological string on $\IC^3$ with a framed brane on one leg of the toric diagram has been studied extensively in the physics literature. For instance, it was first conjectured in \cite{mv}, and then proved in \cite{llz,ophodge}, that the open amplitudes can be written explicitly in terms of Hodge integrals, paralleling the ELSV formula for Hurwitz numbers.

First, recall that we defined in the previous section the open string amplitude with $h$ holes\footnote{Note that there is no dependence on $t$ anymore since $\IC^3$ has no K\"ahler parameter.}
\be
F_{g,h}(x_1, \ldots, x_h;f)=\sum_{w_i} F_{g,w}(f)  x_1^{w_1} \cdots x_h^{w_h},
\ee
where we wrote explicitly the dependence on the framing $f \in
\IZ$ of the brane. These amplitudes can be rewritten as integrals
of differential forms,
\be
F_{g,h}(x_1, \ldots, x_h;f)= \int W_{g}(x_1,\ldots,x_h;f) \rd x_1 \cdots \rd x_h,
\ee
where we introduced, in terms of partitions $\mu$ (in the following we use the notation introduced in section 2)
\be W_{g}(x_1,\dots,x_h;f) = \sum_{\mu | \ell(\mu) = h} z_\mu
W_{g,\mu} (f) m_{\mu}(x). \ee

The amplitudes $W_{g,\mu}(f)$ can be written in terms of Hodge integrals as follows \cite{mv,llz}:
\begin{multline}
\label{mvfor}
W_{g, \mu} (f)= {(-1)^{g + \ell(\mu)} \over |{\rm Aut}(\mu)|}
(f(f+1))^{\ell(\mu)-1}\prod_{i=1}^{\ell(\mu)} \left( { \mu_i \prod_{j=1}^{\mu_i-1} (\mu_i f+ j) \over
\mu_i!} \right) \\
\times \int_{\overline{M}_{g,h}}
{\Lambda^\vee_g (1) \Lambda^\vee_g (-f-1) \Lambda_g^\vee (f)  \over \prod_{i=1}^{\ell(\mu)} (1-\mu_i \psi_i)}.
\end{multline}

What we would like to do now is to use the formalism proposed in \cite{mm,bkmp} to construct recursively the open amplitudes, hence providing an infinite set of relations between the above Hodge integrals. We will then explain how these relations motivate our conjecture \ref{conj1} for Hurwitz numbers.

Recall that the framed mirror curve \eqref{framedc3} reads
\be
-y^{f+1} + y^f - x = 0.
\label{framedc32}
\ee
The open string modulus is given by $x$, hence $y$ is a good local coordinate. In order to apply the recursion relation \eqref{recursion2}, we will need to find $y=y(x)$. We can invert \eqref{framedc32} near $(x,y)=(0,1)$ to get
\be
\label{yx}
y(x) = 1 - \sum_{n=1}^{\infty} { x^n \prod_{j=0}^{n-2} (n f+j) \over n!} .
\ee
We now define the one-forms
\be
\zeta_n(y,f) =\rd y {(1+f) y-f \over y(y-1)} \Bigl( {y(y-1) \over (1+f) y-f}
 {\rd \over \rd y}\Bigr)^{n+1} {1\over (1+f) ((1+f)y-f)}, \qquad n\ge 0.
\ee
For example,
\be
\zeta_0(y,f) =-{\rd y \over (y+f(y-1))^2}, \quad \zeta_1(y,f)= {-f+y(y-2) +fy^2 \over
(y+f(y-1))^4} \rd y.
\ee
It is straightforward to check that we have the following expansion in powers of $x$:
\be
\label{zeta2}
\zeta_n(x,f)=\sum_{\mu=1}^{\infty}{\mu^{n+2} \prod_{j=1}^{\mu-1} (\mu f+ j) \over
\mu!}x^{\mu-1} \rd x.
\ee
It then follows that
\begin{multline}
\label{exprwg}
 W_{g}(x_1, \dots, x_h) \rd x_1 \cdots \rd x_h = (-1)^{g+h} (f(f+1))^{h-1}\\
\times \sum_{n_i=0}^{3g-3+h}   \langle \tau_{n_1} \cdots \tau_{n_h} \Lambda^\vee_g (1) \Lambda^\vee_g (-f-1) \Lambda^\vee_g (f) \rangle
\prod_{i=1}^h \zeta_{n_i}(y_i,f).
\end{multline}

Now, it was proposed in \cite{mm,bkmp} that these differentials can be computed as follows. Take the framed curve \eqref{framedc32}, which gives the expression \eqref{yx} for $y(x)$. Then,
\be
W_0(x) \rd x = \log y(x) {\rd x \over x}, \qquad W_0 (x_1, x_2) = B(y_1,y_2) - {\rd x_1 \rd x_2 \over (x_1 - x_2)^2},
\ee
where the $y_i$ are defined implicitly in terms of $x_i$ through the series $y(x_i)$ given by \eqref{yx}.
The remaining differentials correspond precisely to the differentials generated through the modified recursion relation \eqref{recursion2} with the one-form \eqref{oneform}.

Let us now be a little more explicit. To perform the recursive calculations, we first need to determine the ramification points of the framed curve \eqref{framedc32}. There is only one ramification point of the $x$-projection, which we denote by $\nu$. It is given by
\be
y(\nu) = {f \over f+1}.
\ee
We then find the two points $q$ and $\bar q$ on the curve such that $x(q) = x(\bar q)$ near the ramification point. Write first
\be
y(q) = {f \over f+1} + z, \qquad y(\bar q) = {f \over f+1} + P(z),
\ee
with
\be
P(z) = - z + \CO(z^2).
\ee
Then $\bar q$ must satisfy
\be
-y(q)^{f+1} + y(q)^f = -y({\bar q})^{f+1} + y({\bar q})^f,
\ee
which we can solve to extract $P(z)$ as a power series in $z$. We obtain:
\begin{multline}
\label{pz}
P(z) = -z - \frac{2\,\left( -1 + f^2 \right)}{3\,f} z^2 - \frac{4\,{\left( -1 + f^2 \right) }^2}{9\,f^2} z^3 \\
- \frac{2\,{\left( 1 + f \right) }^3\,\left( -22 + 57\,f - 57\,f^2
+ 22\,f^3 \right)}{135\,f^3} z^4+\cdots
\end{multline}
Using these ingredients, the one-form $\omega(q)$ defined in \eqref{oneform} becomes, in terms of the coordinate $z$:
\be
\label{onec3}
\omega(z) = \left( \log \left({f \over f+1} + z \right) - \log \left({f \over f+1} + P(z) \right) \right) {(f+1)^3 z \rd z \over (f + (f+1) z) ( - 1 + (f+1) z) }.
\ee
The other one-form that we need, $\rd E_q(p)$, was defined, for genus $0$ curves, in \eqref{deq}. Using the results above we obtain, in terms of $z$ and $y$:
\be
\label{dEc3}
\rd E_z (y) = {1 \over 2} \rd y \left( {1 \over y - {f \over f+1} - z} - {1 \over y - {f \over f+1} - P(z)} \right).
\ee
We can now compute explicitly the differentials generated through the recursion \eqref{recursion2}, which becomes
\be
\ba
\label{recursionc3}
W_g(y, y_1 \dots, y_h) &= {\rm Res}_{z=0}{\rd E_{z}(y) \over \omega(z)} \Big( W_{g-1} \left({f \over f+1} + z, {f \over f+1} + P(z), y_1, \dots, y_h\right)\\
&\quad +\sum_{l=0}^g \sum_{J \subset H} W_{g-l}\left({f \over f+1} + z, y_J\right) W_l\left({f \over f+1} + P(z), y_{H\backslash J}\right)\Big).
\ea
\ee
We perform the residue calculation, and express our results in terms of the one-forms $\zeta_n (y,f) \rd y$ introduced earlier. This makes it straightforward to extract the results for the Hodge integrals \eqref{mvfor}, by using the expression \eqref{exprwg} for the differentials $W_g$.

Here are some of the results we obtained.
\be
\ba
\label{results2}
W_0(y_1, y_2, y_3)&= - (f (f+1))^2 \prod_{i=1}^3   \zeta_0(y_i), \\
W_0(y_1, \ldots, y_4)&= (f (f+1))^3 \sum_{i=1}^4 \zeta_1(y_i)  \prod_{j\not= i} \zeta_0(y_j) , \\
W_1(y) & = {1 \over 24} \Big( (1+f+f^2) \zeta_0(y,f) - f(1+f) \zeta_1 (y,f) \Big) , \\
W_1(y_1, y_2)&={1 \over 24}{f (f+1)}\Bigl( -(1+f+f^2) \zeta_0(y_1) \zeta_1(y_2)+ f (1+f)\zeta_0(y_1) \zeta_2(y_2) \\ &
 \qquad \qquad \qquad +(y_1 \leftrightarrow y_2)
+f(1+f)\zeta_1(y_1) \zeta_1(y_2) \Bigr),\\
W_2(y) &= {1 \over 5760} \Bigl( 2 f(f+1) \zeta_1(y,f) -{7} (1+f+f^2)^2 \zeta_2(y,f) \\
& \qquad +12 f (1+2f+2f^2 + f^3) \zeta_3(y,f)-{5} f^2 (1+f)^2 \zeta_4(y,f) \Bigr).\\
\ea
\ee

\subsection{Back to Hurwitz theory}

We are now ready to motivate our Conjecture \ref{conj1} in Hurwitz
theory. It is known that Hurwitz theory can be seen as a
particular limit of topological string theory on the framed vertex
studied in the previous section, where we send the framing $f$ to
infinity. This can be seen in many ways. A simple way is just to
look at the explicit expression for the topological vertex and
take the limit directly \cite{caporaso}. Alternatively, one can
look at the Hodge integral formula for the framed vertex
conjectured in \cite{mv} and proved in \cite{llz,ophodge}. The
(appropriately scaled) limit of infinite framing involves only
Hodge integrals with one Hodge class insertion, and using the ELSV
formula one relates it to Hurwitz numbers.

Here we use this last approach and look at the Hodge integrals expression for the differentials that we computed.
Recall that for the framed vertex, the objects we studied and generated recursively were the differentials \eqref{exprwg}, which we reproduce here for clarity:
\begin{multline}
\label{exprwg2}
 W_{g}(x_1, \dots, x_h) \rd x_1 \cdots \rd x_h = (-1)^{g+h} (f(f+1))^{h-1}\\
\times \sum_{n_i=0}^{3g-3+h}   \langle \tau_{n_1} \cdots \tau_{n_h} \Lambda^\vee_g (1) \Lambda^\vee_g (-f-1) \Lambda^\vee_g (f) \rangle
\prod_{i=1}^h \zeta_{n_i}(x_i,f).
\end{multline}
On the Hurwitz theory side, we computed similar differentials, which were expressed in terms of Hodge integrals in \eqref{elsvgen} using the ELSV formula:
\be
\label{elsvgen2}
 H_{g}(x_1, \dots, x_h) \rd x_1 \cdots \rd x_h  \\
=
\sum_{n_i=0}^{3g-3+h} \langle \tau_{n_1} \cdots \tau_{n_h} \Lambda_g^\vee (1)\rangle
\prod_{i=1}^h \zeta_{n_i}(x_i) ,
\ee
First, using the Mumford relation
\be
\Lambda_g^\vee (t) \Lambda_g^\vee(-t) = (-1)^g t^{2 g},
\ee
we obtain the following relation between the Hodge integrals
\be
\underset{f \to \infty}{\rm lim} {1\over f^{2g}} \langle \tau_{n_1} \cdots \tau_{n_h} \Lambda^\vee_g (1) \Lambda^\vee_g (-f-1) \Lambda^\vee_g (f) \rangle = (-1)^g
\langle \tau_{n_1} \cdots \tau_{n_h} \Lambda_g^\vee (1)\rangle .
\ee
Moreover, using the expansions \eqref{zeta1} and \eqref{zeta2} for the one-forms $\zeta_n(x)$ and $\zeta_n(x,f)$, it is easy to see that
\be
\underset{f \to \infty}{\rm lim} \zeta_n\left( {x \over f},f \right)   = \zeta_n(x).
\ee
Notice that we renormalized the variable $x$ by $x \mapsto x/f$ before taking the limit $f \to \infty$.

As a result, we obtain that
\be
\underset{f \to \infty}{\rm lim} \left( {(-1)^{h} \over f^{2g + 2h - 2} } W_g \left({x_1 \over f}, \ldots, {x_h \over f}\right)  \right)  = H_g(x_1, \ldots, x_h),
\ee
which relates the Gromov-Witten potentials to the Hurwitz generating functions.
This relation can also be seen directly by computing the limits of the expressions \eqref{results2} obtained through the recursion relation.

Finally, as we have seen in the previous subsection, following the proposal of \cite{mm, bkmp} the differentials $W_g(x_1, \ldots, x_h) \rd x_1 \cdots \rd x_h$ can be computed in the recursive formalism discussed previously. This implies that the differentials  $H_g (x_1, \ldots, x_h) \rd x_1 \cdots \rd x_h $ in Hurwitz theory can also be computed through this recursive formalism, after taking the limit $f \to \infty$ appropriately. We can work out this limit explicitly for all the objects entering in the recursive relation. For
example, the curve \eqref{framedc32}, after setting $x\rightarrow x/f$ and $y=1-T/f$, reads,
\be
x=T \Bigl( 1-{T\over f}\Bigr)^f,
\ee
which as $f\rightarrow \infty$ becomes precisely the definition of
the tree function $x=T(x)\re^{-T(x)}$. By working out the infinite
framing limit of the rest of the ingredients, we obtain the
statement of Conjecture \ref{conj1}.

\bibliographystyle{amsalpha}

\end{document}